\documentclass[final]{siamltex}
\usepackage{enumerate}
\usepackage{amsmath, amssymb}
\usepackage{amsfonts}
\usepackage{graphicx}
\usepackage[mathscr]{eucal}
\usepackage{stmaryrd}
\usepackage[ruled,oldcommands]{algorithm2e}
\usepackage{color}
\usepackage{mathtools}
\graphicspath{{},{figures/},{../figures/}}
\DeclareGraphicsExtensions{.pdf}

\newcommand{\mbb}[1]{\mathbb{#1}}

\newcommand{\mbR}{\mathbb{R}}

\newcommand{\scrC}{\mathscr{C}}
\newcommand{\scrF}{\mathscr{F}}
\newcommand{\scrI}{\mathscr{I}}
\newcommand{\scrN}{\mathscr{N}}
\newcommand{\mO}{\mathcal{O}}

\newcommand{\bv}{\mathbf{v}}
\newcommand{\bZ}{\mathbf{Z}}
\newcommand{\bU}{\mathbf{U}}
\newcommand{\bA}{\mathbf{A}}
\newcommand{\bW}{\mathbf{W}}
\newcommand{\bP}{\mathbf{P}}
\newcommand{\bD}{\mathbf{D}}
\newcommand{\bC}{\mathbf{C}}
\newcommand{\bB}{\mathbf{B}}

\newcommand{\bL}{\mathbf{L}}
\newcommand{\bI}{\mathbf{I}}
\newcommand{\bF}{\mathbf{F}}
\newcommand{\bbR}{\hat{\mathbf{R}}}

\newcommand{\bG}{\mathbf{G}}

\newcommand{\bbB}{\mathbf{B}}
\newcommand{\bX}{\mathbf{X}}
\newcommand{\bH}{\mathbf{H}}

\newcommand{\ba}{\mathbf{a}}

\newcommand{\bx}{\mathbf{x}}
\newcommand{\by}{\mathbf{y}}

\newcommand{\bb}{\mathbf{b}}
\newcommand{\bu}{\mathbf{u}}
\newcommand{\bw}{\mathbf{w}}

\newcommand{\bzero}{\mathbf{0}}

\newcommand{\euZ}{\boldsymbol{\mathscr{Z}}}
\newcommand{\euX}{\boldsymbol{\mathscr{X}}}
\newcommand{\euT}{\boldsymbol{\mathscr{T}}}
\newcommand{\bUps}{\mathbf{\Upsilon}}
\newcommand{\bGam}{\mathbf{\Gamma}}
\newcommand{\bPhi}{\mathbf{\Phi}}

\newcommand{\bAhat}{{}\hat{\mathbf{A}}}
\newcommand{\bAtilde}{{}\tilde{\mathbf{A}}}
\newcommand{\bZhat}{{}\hat{\mathbf{Z}}}
\newcommand{\bPhat}{{}\hat{\mathbf{P}}}
\newcommand{\euZhat}{\hat{\boldsymbol{\mathscr{Z}}}}

\title{An adaptive algebraic multigrid algorithm for low-rank canonical tensor decomposition\thanks{This work was supported by NSERC of Canada and was performed in part under the auspices of the U.S. Department of Energy by Lawrence Livermore National Laboratory under Contract DE-AC52-07NA27344.}}

\author{ Hans De Sterck\footnotemark[2]
\and Killian Miller\footnotemark[2]
}

\begin{document}
\maketitle

\renewcommand{\thefootnote}{\fnsymbol{footnote}}
\footnotetext[2]{Department of Applied Mathematics, University of Waterloo,
Waterloo, Ontario, N2L 3G1, Canada (hdesterck@uwaterloo.ca, k7miller@uwaterloo.ca).}
\renewcommand{\thefootnote}{\arabic{footnote}}

\begin{abstract}
This paper presents a multigrid algorithm for the computation of the rank-$R$ canonical decomposition of a tensor for low rank $R$. Standard alternating least squares (ALS) is used as the relaxation method. Transfer operators and coarse-level tensors are constructed in an adaptive setup phase based on multiplicative correction and on Bootstrap algebraic multigrid. An accurate solution is then computed by an additive solve phase based on the Full Approximation Scheme. Numerical tests show that for certain test problems the multilevel method significantly outperforms standalone ALS when a high level of accuracy is required.
\end{abstract}



\pagestyle{myheadings}
\thispagestyle{plain}
\markboth{HANS DE STERCK AND KILLIAN MILLER}{MULTIGRID FOR CANONICAL TENSOR DECOMPOSITION}

\section{Introduction}
\label{sec:intro}

In this paper we present a multigrid method for accurately computing a low-rank canonical decomposition of a tensor. 
An $N$th-order tensor is an $N$-dimensional array of size $I_1 \times \cdots \times I_N$ \cite{TensorReview}. The {\em order} of a tensor is the number of modes (dimensions), and the size of the $n$th mode is $I_n$ for $n = 1,\ldots,N$. The canonical tensor decomposition is a higher-order generalization of the matrix singular value decomposition (SVD) in that it decomposes a tensor as a sum of rank-one components. For example, let $\euT$ be 
an arbitrary $N$th-order tensor of rank $R$, meaning that it can be expressed as a sum of no fewer than $R$ rank-one components. Then its canonical decomposition is given by
\begin{align}\label{eq:tensorCP}
\euT = \sum_{r = 1}^R \ba_r^{(1)}\circ\cdots\circ\ba_r^{(N)},
\end{align}
where $\circ$ denotes the vector outer product. The $r$th rank-one component is formed by taking the vector outer product of $N$ column vectors $\ba_r^{(n)} \in \mbR^{I_n}$ for $n = 1,\ldots,N$. For each mode $n = 1,\ldots,N$, one can store the vectors $\ba_1^{(n)},\ldots,\ba_R^{(n)}$ as the columns of an $I_n \times R$ matrix $\bA^{(n)}$. The matrix $\bA^{(n)}$ is referred to as the mode-$n$ {\em factor matrix} and its columns are the mode-$n$ {\em factors}. The canonical tensor decomposition may then be expressed in terms of the factor matrices by the double bracket notation $\llbracket \bA^{(1)},\ldots,\bA^{(N)} \rrbracket$, which is defined as the summation in \eqref{eq:tensorCP}.
 
We refer to the canonical decomposition as CANDECOMP/PARAFAC (CP) after the names originally given to it in early papers on the subject \cite{CarrolALS,HarshmanALS}. Whereas \eqref{eq:tensorCP} is an example of an exact CP decomposition, referred to as the {\em rank decomposition} since $R = \rank(\euT)$, our goal is to find a CP decomposition for an arbitrary $N$th-order tensor $\euZ$ and a given number of components $R$. Depending on $R$ the CP decomposition may only approximate $\euZ$, in which case it is natural to look for the ``best'' such approximation in some sense.
For example, this problem can be made more concrete by posing it as an optimization problem:
\begin{align}\label{eq:CPfunctional}
\text{minimize} & \quad f(\bA^{(1)},\ldots,\bA^{(N)}) \coloneqq \frac{1}{2}\Big\|\euZ - \llbracket \bA^{(1)},\ldots,\bA^{(N)} \rrbracket \Big\|^2.
\end{align}
In words, we seek the factor matrices that minimize the functional $f$.
Here, $\|\cdot\|$ denotes the Frobenius norm of a tensor, defined as the square root of the sum of the squared entries of all tensor elements. In what follows, for matrices $\|\cdot\|$ refers to the Frobenius norm, and for vectors it refers to the vector two-norm. A general approach to solving the optimization problem \eqref{eq:CPfunctional} is to find a set of nontrivial factor matrices that zero out the gradient of $f$. In other words, at any local minimum of \eqref{eq:CPfunctional} the {\em first-order optimality equations} must be satisfied:
\begin{align}\label{eq:gradf}
\frac{\partial f}{\partial \bA^{(n)}} = \bzero \quad \text{for~} n = 1,\ldots,N.
\end{align}
In \cite{CPOPT}, Acar, Dunlavy and Kolda propose solving the first-order optimality equations by applying gradient-based optimization methods such as the nonlinear conjugate gradient method together with a line search. Along the same vein, a new method for computing the CP decomposition was recently developed by combining a nonlinear generalized minimal residual method with a line search \cite{HansNGMRES}.  
In this paper we follow suit by solving the first-order optimality equations using a multigrid approach.

The optimization problem \eqref{eq:CPfunctional} is non-convex, consequently, it may admit multiple local minima.
Moreover, for any local minimizer there is a continuous manifold of equivalent minimizers \cite{TensorReview}. This manifold arises because of a scaling indeterminacy inherent to the CP decomposition, i.e., the individual factors composing each rank-one term can be rescaled without changing the rank-one term. The CP decomposition also exhibits a permutation indeterminacy in that the rank-one component tensors can be reordered arbitrarily \cite{TensorReview}. In the following sections we discuss how the scaling and permutation indeterminacies can be removed by imposing a specific normalization and ordering of the factors. However, the CP decomposition may still exhibit multiple local minima for some problems, and depending on the initial guess, iterative methods may converge to different stationary points. Moreover, for certain tensors and certain values of $R$ a best rank-$R$ approximation does not exist \cite{Lim}. Despite these difficulties, the exact CP decomposition has been shown to be unique up to scaling and permutation indeterminacies under mild conditions relating the ranks of the factor matrices with the tensor rank, and CP decompositions are used in many application fields \cite{TensorReview}.

The primary application of the CP decomposition is as a tool for data analysis, where it has been used in a variety of fields including chemometrics, data mining, image compression, neuroscience and telecommunications. A second class of problems is related to the decomposition of tensors arising from PDE discretization on high-dimensional regular lattices \cite{Oseledets2006,Oseledets2008,Khoromskij2008,Khoromskij2009}. Many algorithms have been proposed for computing the CP decomposition \cite{TensorReview,Lathauwer2004,Lathauwer2006,Tomasi2006,CPOPT}. However, the workhorse algorithm today is still the original alternating least squares (ALS) method, which was first proposed in 1970 in early papers on the CP decomposition \cite{CarrolALS,HarshmanALS}. The alternating least squares method is simple to implement, and often performs adequately, however, it can be very slow to converge, and its convergence may depend strongly on the initial guess. Despite its simplicity and potential drawbacks, it has proved difficult over the years to develop alternatives to ALS that significantly improve upon it in a robust way for large classes of problems. As a result, ALS-type algorithms are still considered the method of choice in practice.


The multigrid method described in this paper is intended to solve the first-order optimality equations for the CP optimization problem \eqref{eq:gradf}. It consists of two multilevel phases: a multiplicative correction scheme as the setup phase, and an additive correction scheme as the solve phase. In the setup phase a multiplicative correction scheme is used in conjunction with Bootstrap algebraic multigrid (BAMG) interpolation \cite{AchiBAMG,KahlThesis} to not only build the necessary transfer and coarse-level operators, but also to compute initial approximations of the factor matrices. This phase uses the ALS method as the relaxation scheme on all levels. The setup phase is adaptive in the sense that the transfer operators are continually improved using the most recent approximation to the solution factor matrices. In order for the exact solution to be a fixed point of the multiplicative correction scheme, it needs to lie exactly in the range of interpolation at convergence. However, since each interpolation operator attempts to fit multiple factors (in a least-squares sense) this condition can be met only approximately. Therefore, after a few setup cycles we freeze the operators and use additive correction cycles in the solve phase, which can still converge when the exact solution lies only approximately in the range of interpolation.
The combination of a multiplicative setup scheme and BAMG has already been considered in \cite{Kushnir2010}, where it formed the basis of an efficient eigensolver for multiclass spectral clustering problems. A similar approach was also proposed in \cite{KahlThesis,KahlBAMG}. In the solve phase we use the Full Approximation Scheme (FAS) \cite{AchiFAS,McCormick1983,mgtut,Trottenberg} to efficiently obtain an accurate solution.
Our multigrid framework is closely related to recent work on an adaptive algebraic multigrid solver for extremal singular triplets and eigenpairs of matrices \cite{SVDAMG}, and to a lesser degree to multigrid methods for Markov chains \cite{OTF,NSMC,bootstrapMarkov}. Note that in this paper we give a complete but concise description of the proposed algorithm. Readers who are not familiar with adaptive algebraic multigrid (AMG) methods can find more background information about the general adaptive AMG approach in, for example, \cite{BrezinaASA,BrezinaAAM,SVDAMG,OTF}.


The multigrid method proposed in this paper is expected to work well for tensors that have properties which make a multilevel approach beneficial, but not for generic tensors that lack these properties. Just as in the case of multigrid for matrix systems derived from PDE discretizations, our multilevel approach can lead to significant speedup when error components that are damped only weakly by the fine-level process can be represented and damped efficiently on coarser levels. We expect this to be the case for the decomposition of certain higher-order tensors that arise in the context of PDE discretization on high-dimensional regular lattices \cite{Oseledets2006,Oseledets2008,Khoromskij2008,Khoromskij2009}, and we will illustrate the potential benefits of the proposed multigrid method for these types of problems. It should also be noted that the type of multigrid acceleration proposed in this paper will only be effective for low-rank decompositions with small $R$ (e.g., up to 5 or 6). This is because a single interpolation operator is associated with an entire factor matrix, and each interpolation operator can only be expected to represent a small number of factors in a sufficiently accurate way, especially if the desired factors have little in common. These restrictions are entirely analogous to the case of the adaptive multigrid method for computing SVD triplets of a matrix \cite{SVDAMG}.

The remainder of this paper is structured as follows. In \S\ref{sec:notation} we define the basic notation and definitions used throughout this paper. Section \ref{sec:CP} presents the first-order optimality equations and describes the alternating least squares method. Section \ref{sec:multSetup} describes the multilevel setup phase, and \S\ref{sec:FAS} describes the multilevel solve phase. Implementation details and numerical results are presented in \S\ref{sec:numerics} followed by concluding remarks in \S\ref{sec:conc}.

\section{Notation and definitions}
\label{sec:notation}

This section outlines the notation that is used throughout this paper, much of which has been adopted from \cite{CPOPT}. We also review some basic definitions and identities that are important to this paper. For further details we refer to the survey paper by Kolda and Bader \cite{TensorReview}, and the extensive references therein.

Vectors (tensors of order one) are denoted by boldface lowercase letters, e.g., $\bv$. Matrices (tensors of order two) are denoted by boldface capital letters, e.g., $\bA$. Higher-order tensors are denoted by boldface Euler script letters, e.g., $\euZ$. The $i$th entry of a vector $\bv$ is denoted by $v_i$, element $(i,j)$ of a matrix $\bA$ is denoted by $a_{ij}$, and, for example, element $(i,j,k,\ell)$ of a fourth-order tensor $\euZ$ is denoted by $z_{ijk\ell}$. The $j$th column of a matrix $\bA$ is denoted by $\ba_j$. The $n$th element of a sequence is denoted by a superscript in parentheses, e.g., $\bA^{(n)}$. In general, indices range from 1 to their capital versions, e.g., $n = 1,\ldots,N$.


{\em Matricization}, also referred to as unfolding or flattening, is the process of reordering the elements of a tensor into a matrix. In this paper we are only interested in mode-$n$ matricization, which arranges the mode-$n$ {\em fibers} to be the columns of the resulting matrix. Note that a fiber is a higher-order analogue of matrix rows/columns, which is obtained by fixing every index of a tensor but one. Given a tensor $\euZ$ the mode-$n$ matricized version is denoted by $\bZ_{(n)}$.

Matricization provides an elegant way to describe the product of a tensor by a matrix in mode $n$.
The $n$-mode matrix product of a tensor $\euZ \in \mbR^{I_1\times\cdots\times I_N}$ with a matrix $\bA \in \mbR^{J\times I_n}$ is denoted by $\euZ \times_n \bA$ and is of size $I_1 \times \cdots \times I_{n-1}\times J \times I_{n+1}\times\cdots\times I_N$. This product can be expressed in terms of unfolded tensors as follows
\[
\euX = \euZ \times_n \bA \quad\Leftrightarrow\quad \bX_{(n)} = \bA\bZ_{(n)}.
\]

For matrices $\bA \in \mbR^{I\times K}$ and $\bbB \in \mbR^{J\times K}$, their {\em Khatri-Rao product} results in a matrix of size $(IJ) \times K$ defined by
\[
\bA \odot \bbB = [\ba_1 \otimes \bb_1 \cdots \ba_K \otimes \bb_K],
\]
where $\otimes$ is the {\em Kronecker product}. These products have many useful properties, however, we will only need the associativity of the Khatri-Rao product, and the {\em mixed-product property} of the Kronecker product, i.e., $(\bA \otimes \bbB)(\bC \otimes \bD) = \bA\bC \otimes \bbB\bD$. We can then easily prove the following useful identity for any sequences of matrices $\bA^{(n)}$ and $\bbB^{(n)}$, $n = 1,\ldots,N$, of the appropriate sizes
\begin{align}\label{eq:propKKR}
\bA^{(1)}\bbB^{(1)} \odot \cdots \odot \bA^{(N)}\bbB^{(N)} = \left(\bA^{(1)} \otimes \cdots \otimes \bA^{(N)}\right)\left(\bbB^{(1)} \odot \cdots \odot \bbB^{(N)}\right).
\end{align}
Another useful relationship between tensors and their matricized versions is as follows. Let $\euZ \in \mbR^{I_1\times\cdots\times I_N}$ and $\bA^{(n)} \in \mbR^{J_n \times I_n}$ for $n = 1,\ldots,N$. Then, for any $n \in \{1,\ldots,N\}$
\begin{align}\label{eq:propMatricized}
\euX = \euZ \times_1 \bA^{(1)} \cdots \times_N\bA^{(N)} \quad \Leftrightarrow \quad \bX_{(n)} = \bA^{(n)}\bZ_{(n)}\left(\bA^{(N)} \otimes \cdots \otimes \bA^{(1)}\right)^T.
\end{align}
A proof of this property is given in \cite{KoldaMoptsHOD}. 

To denote the product of a tensor and a sequence of matrices over some nonempty subset of the modes $\scrN = \{n_1,\ldots,n_k\} \subset \{1,\ldots,N\}$, we use the following notation as shorthand
\[
\euZ \times_{n \in \scrN} \bA^{(n)} = \euZ \times_{n_1} \bA^{(n_1)} \cdots \times_{n_k}\bA^{(n_k)}.
\]


\section{CP first-order optimality equations and alternating least squares}
\label{sec:CP}


The first-order optimality equations for the CP decomposition are obtained by setting the gradient of the functional in \eqref{eq:CPfunctional} equal to zero. Following \cite{CPOPT}, for each mode $n \in \{1,\ldots,N\}$ the derivative of $f$ with respect to $\bA^{(n)}$ can be written as an $I_n \times R$ matrix
\begin{align}\label{eq:gradEqn}
\bG^{(n)} = -\bZ_{(n)}\bPhi^{(n)} + \bA^{(n)}\bGam^{(n)},
\end{align}
where
\begin{align}\label{eq:Phi}
\bPhi^{(n)} = \bA^{(N)} \odot\cdots\odot\bA^{(n+1)}\odot\bA^{(n-1)}\odot\cdots\odot\bA^{(1)},
\end{align}
and
\begin{align}\label{eq:Gamma}
\bGam^{(n)} = \bUps^{(1)}\ast\cdots\ast\bUps^{(n-1)}\ast\bUps^{(n+1)}\ast\cdots\ast\bUps^{(N)}
\end{align}
with $\bUps^{(n)} = {\bA^{(n)}}^T\bA^{(n)}$ for $n = 1,\ldots,N$. We note that since the $R\times R$ matrix $\bGam^{(n)}$ is the Hadamard (element-wise) product ``$\ast$" of symmetric positive semidefinite matrices, it too must be symmetric positive semidefinite \cite{Styan}. Moreover, if each $\bA^{(n)}$ has full rank then $\bGam^{(n)}$ will be symmetric positive definite (SPD). The first-order optimality equations are then given by
\begin{align}\label{eq:firstorderopt}
\bG^{(n)} = \bzero, \quad n = 1,\ldots,N.
\end{align}

These equations offer a simple way to describe the ALS method for the CP decomposition. One iteration of ALS is equivalent to applying one iteration of block nonlinear Gauss--Seidel (BNGS) to the optimality equations \eqref{eq:firstorderopt}. Iterating through the modes sequentially, at the $n$th step the factor matrices are fixed for all modes except $n$, and the resulting linear least-squares problem is solved for $\bA^{(n)}$. In particular, $\bGam^{(n)}$ and $\bPhi^{(n)}$ are updated and $\bA^{(n)} \leftarrow \bZ_{(n)}\bPhi^{(n)}(\bGam^{(n)})^\dagger$, where $(\bGam^{(n)})^\dagger$ is the Moore--Penrose pseudoinverse of $\bGam^{(n)}$. Owing to the scaling indeterminacy inherent in the CP decomposition, it is possible that during ALS some factors may tend to infinity while others may compensate by tending to zero, such that the rank-one components remain bounded. This behavior can be avoided by using a normalization strategy. After each complete ALS iteration the factors of the $r$th component are normalized according to
\begin{align}\label{eq:ALSnorm}
\ba_r^{(n)} \mapsto \lambda_r\left(\frac{\ba_r^{(n)}}{\|\ba_r^{(n)}\|}\right) \quad \text{for~} n = 1,\ldots,N, \quad \lambda_r = \left(\|\ba_r^{(1)}\|\ldots\|\ba_r^{(N)}\|\right)^{1/N}
\end{align}
for $r = 1,\ldots,R$. This normalization equilibrates the norms of the factors of each component, i.e., $\|\ba_r^{(1)}\| = \cdots = \|\ba_r^{(N)}\|$. The ALS algorithm described here is used as the relaxation method and coarsest-level solver in the setup phase (see \S\ref{sec:multSetup}). 
We note that upon completion of the ALS iterations the rank-one terms are sorted in decreasing order of the normalization factors $\lambda_r$.

\section{Multiplicative setup phase}
\label{sec:multSetup}

In this section we describe the multilevel hierarchy constructed in the setup phase of our solver. We use two-level notation to describe the interaction of two grids at a time. Coarse grid quantities will usually be denoted by a subscript ``$c$'', except in cases where a superscript ``$c$'' improves readability. Fine grid quantities and intergrid transfer operators have neither subscripts nor superscripts. The multiplicative setup phase described here is similar in concept to the setup phase of \cite{SVDAMG} for computing the SVD of a matrix, and more details about this general concept can be found in \cite{SVDAMG}.

\subsection{Derivation of coarse-level equations}
\label{sec:coarseeq}

The fine-level equations are given by the gradient equations stated in \S\ref{sec:CP}, i.e.,
\begin{align}\label{eq:fineeq1}
\bZ_{(n)}\Big(\bA^{(N)} \odot\cdots\odot\bA^{(n+1)}\odot\bA^{(n-1)}\odot\cdots\odot\bA^{(1)}\Big) = \bA^{(n)}\bGam^{(n)} 
\end{align}
for $n = 1,\ldots,N$. Suppose there exist $N$ full rank operators $\bP^{(n)} \in \mathbb{R}^{I_n \times I_{n,c}},~I_{n,c} < I_n$, such that $\bA^{(n)}$ lies approximately in the range of $\bP^{(n)}$, i.e., for each $n$, $\bA^{(n)} \approx \bP^{(n)}\bA_c^{(n)}$ for some coarse-level variables $\bA_c^{(n)} \in \mathbb{R}^{I_{n,c}\times R}$. We note that since each factor matrix has $R$ columns, it is unlikely that we can achieve equality. 
Then it follows that the solution to \eqref{eq:fineeq1} can be approximated by solving a coarse-level problem
\begin{align}\label{eq:coarseeqq}
& {\bP^{(n)}}^T\bZ_{(n)}\Big(\bP^{(N)}\bA_c^{(N)} \odot\cdots\odot\bP^{(n+1)}\bA_c^{(n+1)}\odot\bP^{(n-1)}\bA_c^{(n-1)}\odot\cdots\odot\bP^{(1)}\bA_c^{(1)}\Big)\nonumber\\
& = \Big({\bP^{(n)}}^T\bP^{(n)}\Big)\bA_c^{(n)}\bGam_c^{(n)} \quad \text{for } n = 1,\ldots,N,
\end{align}
followed by interpolation. Here, $\bGam_c^{(n)}$ is defined as in \eqref{eq:Gamma} with 
\[
\bUps_c^{(n)} = {\bA_c^{(n)}}^T({\bP^{(n)}}^T\bP^{(n)})\bA^{(n)} \quad \text{for } n = 1,\ldots,N.
\]
By property \eqref{eq:propKKR} the left-hand side of \eqref{eq:coarseeqq} can be written as
\begin{align*}
{\bP^{(n)}}^T\bZ_{(n)}\Big(\bP^{(N)}\otimes\cdots\otimes\bP^{(n+1)}\otimes\bP^{(n-1)}\otimes\cdots\otimes\bP^{(1)}\Big)\bPhi^{(n)} \quad\text{for } n = 1,\ldots,N,
\end{align*}
with $\bPhi^{(n)}$ given by \eqref{eq:Phi}.
It then follows by property \eqref{eq:propMatricized} that we can write the coarse-level tensor as
\begin{align}\label{eq:ctensor}
\euZ^c = \euZ \times_1 \bP_1^T \times_2 \bP_2^T \cdots \times_N \bP_N^T,
\end{align}
which is essentially a higher dimensional analogue of the Galerkin coarse-level operator that is commonly used in multigrid for the matrix case. Defining $\bB^{(n)} = {\bP^{(n)}}^T\bP^{(n)}$ for each mode $n$ we can write the coarse-level version of \eqref{eq:fineeq1} as follows
\begin{align}\label{eq:coarseeq1}
\bZ_{(n)}^c\Big(\bA_c^{(N)} \odot\cdots\odot\bA_c^{(n+1)}\odot\bA_c^{(n-1)}\odot\cdots\odot\bA_c^{(1)}\Big) = \bB^{(n)}\bA_c^{(n)}\bGam_c^{(n)}.
\end{align}
By the full rank assumption on the interpolation operators it follows that $\bB^{(n)}$ is SPD, and hence we can compute its Cholesky factorization $\bB^{(n)} = \bL^{(n)}{\bL^{(n)}}^T$, where $\bL^{(n)}$ is an $I_{n,c} \times I_{n,c}$ nonsingular lower triangular matrix. 
The Cholesky factors can be used to transform \eqref{eq:coarseeq1}, whereby one obtains an equivalent set of equations that correspond to the first-order optimality equations of a coarse-level CP minimization problem.
Making the change of variables $\bAhat_c^{(n)} = {\bL^{(n)}}^T\bA_c^{(n)}$ for each mode $n$, and by appealing to property \eqref{eq:propKKR}, it follows that \eqref{eq:coarseeq1} can be written as
\[
\bZhat_{(n)}^c\Big(\bAhat_c^{(N)} \odot\cdots\odot\bAhat_c^{(n+1)}\odot\bAhat_c^{(n-1)}\odot\cdots\odot\bAhat_c^{(1)}\Big) = \bAhat_c^{(n)}\hat{\bGam}_c^{(n)},
\]
where
\[
\bZhat_{(n)}^c = {\bL^{(n)}}^{-1}\bZ_{(n)}^c\Big({\bL^{(N)}}^{-1} \otimes\cdots\otimes{\bL^{(n+1)}}^{-1}\otimes{\bL^{(n-1)}}^{-1}\otimes\cdots\otimes{\bL^{(1)}}^{-1}\Big)^T,
\]
and $\hat{\bGam}_c^{(n)} = \bGam_c^{(n)}$ for $n = 1,\ldots,N$.
Moreover, by property \eqref{eq:propMatricized} the transformed coarse-level tensor is given by
\begin{align}\label{eq:ctensor_transformed}
\euZhat^c & = \euZ \times_1 \Big(\bP^{(1)}{\bL^{(1)}}^{-T}\Big)^T \cdots \times_N \Big(\bP^{(N)}{\bL^{(N)}}^{-T}\Big)^T\nonumber\\
& = \euZ\times_1{\bPhat^{(1)}}^T\cdots \times_N{\bPhat^{(N)}}^T,
\end{align}
where $\bPhat^{(n)} = \bP^{(n)}{\bL^{(n)}}^{-T}$ for $n = 1,\ldots,N$.
Therefore, the coarse-level equations are equivalent to the gradient equations of the following coarse-level functional:
\begin{align}\label{eq:cfunc}
\hat{f}_c(\bAhat_c^{(1)},\ldots,\bAhat_c^{(N)}) \coloneqq \frac{1}{2}\Big\|\euZhat^c - \llbracket \bAhat_c^{(1)},\ldots,\bAhat_c^{(N)} \rrbracket \Big\|^2.
\end{align}
Therefore, the coarse-level equations can be solved by applying ALS to minimize the coarse-level functional $\hat{f}_c$. Initial approximations of the transformed coarse-level variables may be obtained by applying a restriction operator $\bbR^{(n)}$ to the current fine-level approximations.
In this paper we follow the standard approach of defining the restriction operators as the transposes of the interpolation operators, i.e.,
\begin{align}
\bbR^{(n)} = {\bPhat^{(n)}}^T \quad \text{for~} n = 1,\ldots,N.
\end{align}
After solving the coarse-level equations, the coarse-grid-corrected fine-level approximations are obtained via prolongation:
\begin{align}\label{eq:CGC}
{\bA_{\text{\tiny CGC}}^{(n)}} = \bPhat^{(n)}\bAhat_c^{(n)}\quad \text{for } n = 1,\ldots,N.
\end{align}
We note that if $\bP^{(n)}$ contained $\bA^{(n)}$ exactly in its range, and if the coarse-level equations were solved exactly, then the coarse-grid-corrected solution would equal the exact fine-level solution:
\[
{\bA_{\text{\tiny CGC}}^{(n)}} = \bPhat^{(n)}\bAhat_c^{(n)} = \bP^{(n)}{\bL^{(n)}}^{-T}{\bL^{(n)}}^T\bA_c^{(n)} = \bP^{(n)}\bA_c^{(n)} = \bA^{(n)}\quad \text{for } n = 1,\ldots,N.
\]
However, since these conditions only hold approximately, we expect \eqref{eq:CGC} to yield an improved but not exact approximation to the fine-level solution. In particular, since the approximation properties of the interpolation operators deteriorate as the number of components is increased, we expect our method to perform well for a relatively small number of components $R$.

\subsection{Bootstrap AMG V-cycles}
\label{sec:bootstrap}

In this section we describe how we use Bootstrap AMG \cite{AchiBAMG,KahlThesis} to find initial approximations of the desired factor matrices, and adaptively determine the interpolation operators that approximately fit the factor matrices. We follow the approach outlined in \cite{Kushnir2010,KahlThesis,SVDAMG}.

We begin by describing the initial BAMG V-cycle. On the finest level we randomly choose $n_t$ {\em test blocks}, where each test block is a collection of $N$ randomly generated {\em test factor matrices} (TFMs) $\bA_t^{(1)},\ldots,\bA_t^{(N)}$. 
The reason we must consider test blocks instead of simply adding more columns to the factor matrices is that rank-one components of the best rank-$R$ CP tensor approximation must be found simultaneously \cite{TensorReview}; contrary to the best rank-$R$ matrix approximation, one cannot obtain the best rank-$R$ CP approximation by truncating the best rank-$Q$ approximation with $Q > R$. We also start with a collection of $N$ randomly generated {\em boot factor matrices} (BFMs) $\bA_b^{(1)},\ldots,\bA_b^{(N)}$, which serve as our initial guess to the desired factor matrices. 
We note that the subscripts ``$t$'' and ``$b$'' serve only to distinguish between the test and boot factors.

In the downward sweep of the first BAMG V-cycle we relax on each of the test blocks separately and also on the BFMs using the ALS algorithm described in \S\ref{sec:CP}. We then coarsen each mode on the finest level and determine the interpolation operators $\bP^{(1)},\ldots,\bP^{(N)}$. The $n$th interpolation operator $\bP^{(n)}$ fits the factors in the $n$th TFMs across all test blocks in a least-squares sense, such that these factors lie approximately in the range of $\bP^{(n)}$. We also build the coarse-level operator $\hat{\euZ}^c$ and restrict the TFMs and BFMs to the first coarse level. This process is then repeated recursively until some coarsest level is reached, from which point on we relax only on the BFMs.

In the upward sweep of the first cycle, starting from the coarsest level, we recursively interpolate the BFMs up to the next finer level which gives the coarse-grid-corrected approximation on that level. We then relax the coarse-grid correction (CGC) using ALS. This process continues until the CGC on the finest level has been relaxed by ALS.

The initial BAMG V-cycle can be followed by several additional BAMG V-cycles. These cycles are the same as the initial cycle except for one key difference. In the downward sweep the $n$th interpolation operator $\bP^{(n)}$ fits the factors in the $n$th TFMs across all test blocks as well as the factors in the $n$th BFM. Since the BFMs serve as our initial approximation for the additive phase of the algorithm, they must be well represented by interpolation if the additive solve phase is to converge. 

\subsection{Interpolation sparsity structure: coarsening}
\label{sec:coarsening}

Construction of the interpolation operators proceeds in two phases. In the first phase, the sparsity structure of $\bP^{(n)}$ is determined by selecting a subset of the fine-level indices $\Omega_n = \{1,\ldots,I_n\}$. This subset, denoted by $\scrC_n$, is the set of coarse indices for mode $n$ (it has cardinality $I_{n,c} < I_n$). 
Fine-level points that are not selected to the coarse level are represented by the set of fine indices $\scrF_n = \Omega_n\setminus\scrC_n$. 
For each point $i \in \scrF_n$ we define a set of coarse interpolatory points $\scrC_n^i$, which contains coarse points that $i$ interpolates from. For convenience we assume that the points in $\scrC_n^i$ are labeled by their coarse-level indices. Furthermore, for any fine-level point $i \in \scrC_n$ we let $\alpha(i)$ denote its coarse-level index. The interpolation operator $\bP^{(n)}$ is defined by
\[
p_{ij}^{(n)} = \left\{
\begin{array}{ll}
w_{ij}^{(n)}, & i \in \scrF_n \text{ and } j \in \scrC_n^i\\[1mm]
1, & i \in \scrC_n \text{ and } j = \alpha(i)\\[1mm]
0, & \text{otherwise,}
\end{array}\right.
\]
where the $w^{(n)}_{ij}$s are the interpolation weights for mode $n$. The interpolation weights are determined by a least-squares process described in \S\ref{sec:LSQweights}. In this paper we use standard geometric coarsening for each mode, whereby $\scrC_n$ consists of the odd numbered points in $\Omega_n$, and $\scrF_n$ consists of the even numbered points. It follows that $\alpha(i) = (i+1)/2$ in this case. For each $i \in \scrF_n$ we define $\scrC_n^i = \{\alpha(i-1),\alpha(i+1)\}$ (coarse-level labels) except possibly at the right endpoint. This coarsening works well when the modes have approximately the same size, however, for tensors in which the sizes of some modes vary widely, a more aggressive coarsening for the larger modes may be considered. In \S\ref{sec:CPAMGmult} we discuss a straightforward approach to coarsening tensors with varying mode sizes. While the simple coarsening procedure discussed here works well for the test problems considered in \S\ref{sec:numerics} (PDE problems on high-dimensional regular lattices), more general coarsening algorithms for other types of tensors would be desirable. The development of such algorithms remains an interesting topic of future research.


\subsection{Least-squares determination of interpolation weights}
\label{sec:LSQweights}

Suppose that mode $n$ has been coarsened and that $\scrC_n$ and $\scrF_n$ are given. Further suppose that the factors in the $n$th test factor matrices across all test blocks are stored as the columns of the $I_n \times Rn_t$ matrix $\bU_t$, and let $\bU_b = \bA_b^{(n)}$ for the boot factor matrices. Following the approaches of \cite{AchiFAS,Kushnir2010,KahlThesis,SVDAMG} we use a least-squares (LS) process to determine the interpolation weights in the rows of $\bP^{(n)}$ that correspond to points in $\scrF_n$. We want to fit the weights of $\bP^{(n)}$ such that the vectors in $\bU_t$ and $\bU_b$ (except in the first cycle) lie approximately in the range of $\bP^{(n)}$. Let the columns of $\bU_f = [\bU_t~|~\bU_b]$ hold the $n_f = R(n_t+1)$ vectors to be fitted. Let $\bu_k$ be the $k$th column of $\bU_f$. Let $\bu_{k,c}$ be the coarse-level version of $\bu_k$ obtained by injection, and let $(\bu_{k,c})_j$ be its value in the coarse-level point $j$. Also, let $u_{ik}$ be the value of $\bu_k$ in the fine-level point $i$. The interpolation weights of each row that corresponds to a point in $\scrF_n$ may now be determined consecutively by independent least squares fits. Consider a fixed point $i \in \scrF_n$ with coarse interpolatory set $\scrC_n^i$. We solve the following least-squares problem to determine the unknown interpolation weights $w^{(n)}_{ij}$,
\begin{align}\label{eq:LSQinterp}
u_{ik} = \sum_{j \in \scrC_n^i} w^{(n)}_{ij}(\bu_{k,c})_j	\quad\text{for } k = 1,\ldots,n_f.
\end{align}
We make \eqref{eq:LSQinterp} over-determined in all cases by choosing $n_t > M_s/R$, where $M_s$ is the maximum interpolation stencil size for any $i$ on any level, i.e., $|\scrC_n^i| \leq M_s$.
Owing to the standard geometric coarsening of each mode $M_s = 2$, and so it is sufficient to use $n_t = 2$ for any number of components $R > 1$. For a rank-one decomposition we must take $n_t \geq 3$. 

In practice it is necessary to formulate \eqref{eq:LSQinterp} as a weighted least-squares problem, where the weights should bias the fit toward the boot factors. We do so as follows. For a fixed $n$, we compute the gradient $\bG^{(n)}$ according to \eqref{eq:gradEqn} using the factor matrices in the current test block. The weights $w_r$ for the mode-$n$ factor vectors are then defined as
\begin{align}\label{eq:interpweights}
w_r = \frac{\|\bA_t^{(n)}\|^2}{\|\bG^{(n)}\|^2} \quad\text{for } r = 1,\ldots,R.
\end{align}
Note that for a given mode $n$ and test block all factor vectors have the same weights.
The weights for all test blocks are stored in the vector $\bw_t$ of length $Rn_t$. The weights for the boot factors are computed similarly and are stored in the vector $\bw_b$ of length $R$. The full vector of weights $\bw \in \mbR^{n_f}$ is obtained by ``stacking'' $\bw_t$ on top of $\bw_b$. Equation \eqref{eq:interpweights} stems from the observation that $\bG^{(n)}$ is a residual for the $n$th factor matrix. Therefore, since the BFMs should converge much faster than the TFMs, the gradient norm for the BFMs should be much smaller, and hence the weights should be larger. We note that weights corresponding to a single factor matrix are chosen identical in \eqref{eq:interpweights} since we do not want preferential treatment given to different factor vectors, but rather to entire factor matrices. In our implementation the small weighted least-squares problems with diagonal weight matrix $\bW = \diag(\bw)$ are solved by a standard normal equations approach.

\subsection{Pseudocode}
\label{sec:CPAMGmult}

A pseudocode description for a multiplicative setup phase V-cycle with ALS as the relaxation scheme and coarsest-level solver is given by Algorithm \ref{alg:multmlevel}. For the unfamiliar reader we note that a multigrid V-cycle is obtained by performing one recursive solve on each level. If two recursive solves are performed on each level then a multigrid W-cycle is obtained. In this paper we only consider V-cycles for the setup phase.

\begin{algorithm}
\label{alg:multmlevel}
\dontprintsemicolon
\Setnlsty{}{}{.}
\SetNoLine
\caption{V-cycle for setup phase of CP decomposition (CP-AMG-mult)}
\KwIn{tensor $\euZ$, BFMs $\bA^{(1)},\ldots,\bA^{(N)}$, TFMs}
\KwOut{updated BFMs $\bA^{(1)},\ldots,\bA^{(N)}$, updated TFMs}
\BlankLine
\nl Compute the set $\scrI_c = \{n : I_n > I_{n,coarsest}\}$\;
\eIf{$\scrI_c \neq \emptyset$}{
\nl Apply $\nu_1$ relaxations to TFMs in each test block and to $\bA^{(1)},\ldots,\bA^{(N)}$\;
	\For{$n \in \scrI_c$}{
	\nl Build the interpolation operator $\bP^{(n)}$ (on first cycle only use TFMs)\;
	\nl Let $\bB^{(n)} \leftarrow {\bP^{(n)}}^T\bP^{(n)}$\;
	\nl Compute the Cholesky factorization $\bB^{(n)} = \bL^{(n)}{\bL^{(n)}}^T$\;
	\nl Let $\bPhat^{(n)} \leftarrow \bP^{(n)}{\bL^{(n)}}^{-T}$ and $\bbR^{(n)} \leftarrow {\bPhat^{(n)}}^T$\;
	}
\nl Compute the coarse BFMs and coarse TFMs according to \eqref{eq:coarseapprox_new}\;
\nl Compute the coarse-level tensor $\euZhat^c \leftarrow \euZ \times_{n \in \scrI_c} \bbR^{(n)}$\;
\nl Recursive solve:
\[
\{\bAhat_c^{(1)},\ldots,\bAhat_c^{(N)}\} \leftarrow \text{CP-AMG-mult}(\euZhat_c,\bAtilde_c^{(1)},\ldots,\bAtilde_c^{(N)},\text{coarse TFMs})\;
\]
\nl Compute the CGC $\bA^{(n)}$ for $n = 1,\ldots,N$ according to \eqref{eq:CGC_new}\; 
\nl Apply $\nu_2$ relaxations to $\bA^{(1)},\ldots,\bA^{(N)}$\;
}{
\nl Apply $\nu_c$ relaxations to $\bA^{(1)},\ldots,\bA^{(N)}$\;
}
\end{algorithm}

The CP-AMG-mult algorithm recursively coarsens each mode until it reaches some predefined coarsest level. However, since the size of each mode $I_1,\ldots,I_N$ may differ, it is possible that some modes may reach their coarsest level sooner than others. We address this issue as follows. For each mode $n$ we define a threshold $I_{n,coarsest}$ to be the maximum size of that mode's coarsest level. 
For any mode $n$ such that $I_n > I_{n,coarsest}$ we continue coarsening, otherwise we do not coarsen that mode further. Let the modes that still require further coarsening be indexed by the set $\scrI_c = \{n : I_n > I_{n,coarsest}\}$, and let $\scrI_c'$ denote its complement. Then at any given level for each $n \in \scrI_c'$ it follows that $\bPhat^{(n)} = \bI^{(n)}$, where $\bI^{(n)}$ is the $I_n \times I_n$ identity matrix. Setting $\bPhat^{(n)}$ equal to the identity for all $n\in\scrI_c'$ has the following implications. The coarse-level tensor can be computed by taking the product in \eqref{eq:ctensor_transformed} over the modes in $\scrI_c$, instead of for all $n = 1,\ldots,N$. The coarse-level approximations of the BFMs are given by
\begin{align}\label{eq:coarseapprox_new}
\bAtilde_c^{(n)} = \left\{
\begin{array}{ll}
\bbR^{(n)}\bA^{(n)}, & n \in \scrI_c\\[1mm]
\bA^{(n)}, & n \in \scrI_c'
\end{array}\right.
\quad \text{for } n = 1,\ldots,N.
\end{align}
Similarly, the coarse-level approximations of the TFMs in each test block are computed by restricting only those factor matrices indexed by $\scrI_c$. Additionally, the coarse-grid-corrected BFMs are given by
\begin{align}\label{eq:CGC_new}
{\bA_{\text{\tiny CGC}}^{(n)}} = \left\{
\begin{array}{ll}
\bPhat^{(n)}\bAhat_c^{(n)}, & n \in \scrI_c\\[1mm]
\bAhat_c^{(n)}, & n \in \scrI_c'
\end{array}\right.
\quad \text{for } n = 1,\ldots,N.
\end{align}

The size of the coarsest level plays an important role in how the multigrid method performs. If the coarsest level is too large then not enough work is done on the coarser levels and convergence will be slow. Conversely, choosing too small a coarsest level may negatively impact convergence, or in some cases may even cause divergence (as in \cite{Kushnir2010,KahlThesis,SVDAMG}). In practice we find that choosing $I_{n,coarsest} \geq R$ for all $n$ works well.

\section{Full approximation scheme additive solve phase}
\label{sec:FAS}

The Full Approximation Scheme (FAS) \cite{AchiFAS} is the nonlinear analogue of the linear additive correction multigrid method. When applied to linear problems FAS reduces to the usual additive method, and so it is a more general multigrid solver. In this section we discuss how FAS can be used to obtain an additive correction method for the CP decomposition. Two-level notation is used to describe the interaction of two grids at a time with coarse-grid quantities denoted by a subscript ``$c$''.

\subsection{Coarse-level equations}
\label{sec:FAScoarse}

Recall the finest-level equations \eqref{eq:fineeq1}, which can be expressed as
\begin{align*}
\bA^{(n)}\bGam^{(n)} - \bZ_{(n)}\Big(\bA^{(N)} \odot\cdots\odot\bA^{(n+1)}\odot\bA^{(n-1)}\odot\cdots\odot\bA^{(1)}\Big) = \bzero
\end{align*}
for $n = 1,\ldots,N$. Suppose we define nonlinear operators $\bH^{(1)},\ldots,\bH^{(N)}$ such that for any $n \in \{1,\ldots,N\}$
\[
\bH^{(n)} \colon \mbb{R}^{I_1\times R} \times \cdots\times \mbb{R}^{I_N\times R} \rightarrow \mbb{R}^{I_n\times R}
\] 
where
\begin{align}\label{eq:FASop}
\bH^{(n)}(\{\bA\}) \coloneqq \bA^{(n)}\bGam^{(n)}-\bZ_{(n)}\Big(\bA^{(N)} \odot\cdots\odot\bA^{(n+1)}\odot\bA^{(n-1)}\odot\cdots\odot\bA^{(1)}\Big).
\end{align}
Note that we use $\{\bA\}$ as shorthand for $\bA^{(1)},\ldots,\bA^{(N)}$. Then the fine-level problem can be formulated as a system of nonlinear equations
\begin{align}\label{eq:FASfineprob}
&\bH(\{\bA\}) \coloneqq (\bH^{(1)}(\{\bA\}),\ldots,\bH^{(N)}(\{\bA\})) = (\bF^{(1)},\ldots,\bF^{(N)}),
\end{align}
where $\bF^{(n)} = \bzero$ for $n = 1,\ldots,N$ on the finest level. In order to apply the full approximation scheme we require a coarse version of \eqref{eq:FASfineprob}. For each mode $n$ we define the coarse operator
\begin{align}\label{eq:FAScop}
\bH_c^{(n)}(\{\bA_c\}) \coloneqq \bA_c^{(n)}\bGam_c^{(n)}-\bZhat^c_{(n)}\Big(\bA_c^{(N)} \odot\cdots\odot\bA_c^{(n+1)}\odot\bA_c^{(n-1)}\odot\cdots\odot\bA_c^{(1)}\Big),
\end{align}
where $\euZhat_c$ is the coarse-level tensor computed in the multiplicative setup phase. Then the coarse-level FAS equations are given by
\begin{align}\label{eq:FAScoarseprob}
\bH_c(\{\bA_c\}) \coloneqq (\bH_c^{(1)}(\{\bA_c\}),\ldots,\bH_c^{(N)}(\{\bA_c\})) = (\bF_c^{(1)},\ldots,\bF_c^{(N)})
\end{align}
where
\begin{align}\label{eq:FASrhs}
\bF_c^{(n)} = \bbR^{(n)}(\bF^{(n)}-\bH^{(n)}(\{\bA\})) + \bH_c^{(n)}(\{\bAtilde_c\}) \quad \text{for } n = 1,\ldots,N,
\end{align}
and $\bbR^{(n)}$ is the mode $n$ restriction operator from the multiplicative setup phase.
Here $\bAtilde_c^{(n)}$ is the coarse-level approximation of $\bA^{(n)}$ obtained by restriction. Solving \eqref{eq:FAScoarseprob} for $\{\bA_c\}$, the coarse-grid-corrected approximations on the fine level are given by
\begin{align}\label{eq:FASCGC}
{\bA_{\text{\tiny CGC}}^{(n)}} = \bA^{(n)}+\bPhat^{(n)}(\bA_c^{(n)}-\bAtilde_c^{(n)}) \quad \text{for~} n = 1,\dots,N,
\end{align}
where $\bPhat^{(n)}$ is the interpolation operator from the multiplicative setup phase. Together, Equations \eqref{eq:FASfineprob} and \eqref{eq:FAScoarseprob} to \eqref{eq:FASCGC} describe a FAS two-level coarse-grid correction scheme for the CP optimality equations. In the following sections we describe the relaxation scheme used in the solve phase, and give a pseudocode description of the multilevel CP-FAS algorithm.

\subsection{Relaxation}
\label{sec:FASrelax}

We employ block nonlinear Gauss--Seidel (BNGS) as the relaxation scheme and coarsest-level solver for the CP-FAS algorithm (Algorithm \ref{alg:CPFAS}). Applying BNGS to the equations in \eqref{eq:FASfineprob} is similar to applying ALS to the CP optimality equations. One iteration consists of iterating through the modes sequentially, where at the $n$th step $\bGam^{(n)}$ and $\bPhi^{(n)}$ are computed, and then $\bA^{(n)}$ is updated by solving
\begin{align}\label{eq:relaxinnersolve}
\bA^{(n)}\bGam^{(n)} =  \bZ_{(n)}\bPhi^{(n)}+\bF^{(n)}.
\end{align}
When considering how to solve \eqref{eq:relaxinnersolve} for mode $n$, on any level, we recall the following fact: the exact solution is a fixed point of FAS if it is a fixed point of the relaxation scheme and the coarsest-level solver. Suppose we update $\bA^{(n)}$ by post-multiplying the right-hand side of \eqref{eq:relaxinnersolve} by $(\bGam^{(n)})^\dagger$, which is a small $R\times R$ matrix. If $\bGam^{(n)}$ is nonsingular then its pseudoinverse is equivalent to its inverse, in which case there exists a unique solution and the fixed point is preserved. However, if $\bGam^{(n)}$ is singular then post-multiplying by the pseudoinverse will in general not preserve the fixed point. 
Therefore, we consider an alternative approach. 
We propose using a few iterations of Gauss--Seidel (GS) to update $\bA^{(n)}$, which guarantees the fixed point property of our relaxation method. Moreover, a result by Keller \cite{Keller} for positive semidefinite matrices states that if $\bGam^{(n)}$ has nonzero entries on its diagonal then GS must converge to a solution (there may be many) of \eqref{eq:relaxinnersolve}. Owing to the structure of $\bGam^{(n)}$ this condition is equivalent to the fundamental requirement that the factor matrices have nonzero columns. Therefore, we can be confident that GS will converge regardless of whether or not $\bGam^{(n)}$ is singular. In practice we find that only a few GS iterations are necessary to obtain a sufficiently accurate solution to \eqref{eq:relaxinnersolve}, and that further iterations do little to improve the relaxed approximation. In this paper we use ten GS iterations.

After each iteration of BNGS on the finest level the factor matrices are normalized according to \eqref{eq:ALSnorm}. Due to the structure of the FAS equations, in particular the right-hand side in \eqref{eq:FASfineprob}, the scaling indeterminacy is not present on the coarser levels and so normalizing there is unnecessary. We note that the permutation indeterminacy is also removed on the coarser levels because of the (nonzero) right-hand side. Therefore, the rank-one terms are sorted in decreasing order of the normalization factors $\lambda_r$ only on the finest level.


\subsection{CP-FAS algorithm}
\label{sec:CPFAS}

A pseudocode description for an additive solve phase V-cycle is given in Algorithm \ref{alg:CPFAS}. While other cycling schemes such as W-cycles and F-cycles \cite{Trottenberg} are available to us, in this paper we use V-cycles.
We assume that at any given level the current tensor, the index set $\scrI_c$, and the interpolation/restriction operators from the setup phase are available to the algorithm. Note that the parameters $(\nu_1,\,\nu_2,\,\nu_c)$ may be different from those used during the setup phase.\\

\begin{algorithm}[H]
\label{alg:CPFAS}
\dontprintsemicolon
\Setnlsty{}{}{.}
\SetNoLine
\caption{V-cycle for solve phase of CP decomposition (CP-FAS)}
\KwIn{right-hand side matrices $\{\bF\}$, factor matrices $\{\bA\}$}
\KwOut{updated factor matrices $\{\bA\}$}
\BlankLine

\eIf{not on the coarsest level}{
\nl Apply $\nu_1$ relaxations to $\bH(\{\bA\}) = (\bF^{(1)},\ldots,\bF^{(N)})$\;
\nl Coarse initial guess:
\[
\bAtilde_c^{(n)} \leftarrow \left\{
\begin{array}{ll}
\bbR^{(n)}\bA^{(n)}, & n \in \scrI_c\\[1mm]
\bA^{(n)}, & n \in \scrI_c'
\end{array}\right.
\quad n = 1,\ldots,N\;
\]
\nl Coarse right-hand side:
\[
\bF_c^{(n)} \leftarrow \left\{
\begin{array}{ll}
\bH_c^{(n)}(\{\bAtilde_c\})+\bbR^{(n)}(\bF^{(n)}-\bH^{(n)}(\{\bA\})), & n \in \scrI_c\\[1mm]
\bF^{(n)}, & n \in \scrI_c'
\end{array}\right.
\quad\!\! n = 1,\ldots,N\;
\]
\nl Recursive solve:
\[
\{\bA_c^{(1)},\ldots,\bA_c^{(N)}\} \leftarrow \text{CP-FAS}(\{\bF_c\},\{\bAtilde_c\})\;
\]
\nl Coarse-grid correction:
\[
\bA^{(n)} \leftarrow \bA^{(n)} + \left\{
\begin{array}{ll}
\bPhat^{(n)}(\bA_c^{(n)}-\bAtilde_c^{(n)}), & n \in \scrI_c\\[1mm]
\bA_c^{(n)}-\bAtilde_c^{(n)}, & n \in \scrI_c'
\end{array}\right.
\quad n = 1,\ldots,N\;
\]
\nl Apply $\nu_2$ relaxations to $\bH(\{\bA\}) = (\bF^{(1)},\ldots,\bF^{(N)})$\;
}{
\nl Apply $\nu_c$ relaxations to $\bH(\{\bA\}) = (\bF^{(1)},\ldots,\bF^{(N)})$\;
}
\end{algorithm}

It is instructive to mention the differences of this additive solution phase and the additive phase in the SVD case \cite{SVDAMG}. For the SVD, singular vectors can be computed in separate V-cycles and FAS is not required since the singular values are updated in a top-level Ritz step. In the tensor case, all factor vectors need to be computed simultaneously in a single FAS V-cycle, and the weights $\lambda_r$ from \eqref{eq:ALSnorm}, which are in some sense equivalent to the singular values, are updated in these FAS cycles as well, making a Ritz step unnecessary.

\subsection{Full multigrid FAS cycles}
\label{sec:FMG}

For some problems the initial guess provided by the multiplicative setup phase may be too far from the solution to yield a convergent additive solve phase. One way in which we can try to obtain a better initial guess to the fine-level problem is to use Full Multigrid (FMG) \cite{Borzi2005,mgtut,Trottenberg}. Full multigrid is based on {\em nested iterations} whereby coarse levels are used to obtain improved initial guesses for fine-level problems. At any given level the problem is first solved on the next coarser level after which the solution is interpolated to the current level to provide a good initial guess. This process naturally starts at the coarsest level and terminates at the finest. Once an initial guess to the finest-level problem has been obtained we can apply repeated CP-FAS cycles to obtain an improved approximate solution. We use CP-FAS V-cycles as the solver on each level of the FMG cycle, except on the coarsest level where ALS is used (see \S\ref{sec:CP}). A pseudocode description of the FMG-CP-FAS algorithm is given in Algorithm \ref{alg:FMGCPFAS}. We assume that at any given level the current tensor, the index set $\scrI_c$, and the interpolation/restriction operators from the setup phase are available. In Algorithm \ref{alg:FMGCPFAS} we use a subscript $\ell$ to index the current level, where $\ell = 0,\ldots,L-1$. Note that a subscript $\ell$ on an interpolation operator indicates that level $\ell$ is mapped to level $\ell-1$.

\begin{algorithm}[H]
\label{alg:FMGCPFAS}
\dontprintsemicolon
\Setnlsty{}{}{.}
\SetNoLine
\caption{Full multigrid cycle for solve phase of CP decomposition (FMG-CP-FAS)}
\KwOut{finest-level factor matrices $\bA_0^{(1)},\ldots,\bA_0^{(N)}$}
\BlankLine

\nl On the coarsest level apply $\nu$ iterations of ALS with a random initial guess to obtain $\bA_{L-1}^{(1)},\ldots,\bA_{L-1}^{(N)}$\;
\nl Set $\ell \leftarrow L-1$\;

\While{$\ell \neq 0$}{
\nl $\bA_{\ell-1}^{(n)} \leftarrow \bPhat_\ell^{(n)}\bA_\ell^{(n)}$ for $n = 1,\ldots,N$\;
\nl $\{\bA_{\ell-1}^{(1)},\ldots,\bA_{\ell-1}^{(N)}\} \leftarrow \text{CP-FAS}(\{\bzero\},\bA_{\ell-1}^{(1)},\ldots,\bA_{\ell-1}^{(N)})$\;
\nl $\ell \leftarrow \ell-1$\;
}
\end{algorithm}

\section{Implementation details and numerical results}
\label{sec:numerics}

In this section we present the results of numerical tests. All experiments are performed using MATLAB version 7.5.0.342 (R2007b) and version 2.4 of the Tensor Toolbox \cite{TensorToolbox}. Timings are reported for a laptop running Windows XP, with a 2.50 GHz Intel Core 2 Duo processor and 4GB of RAM. Initial guesses for the boot factors and test factors are randomly generated from the standard uniform distribution. The initial boot factors are also used as the initial guess for the standalone ALS method.
The stopping criterion for the numerical tests is based on the gradient of $f$. In particular, with
\begin{align}\label{eq:gradConv}
g(\bA^{(1)},\ldots,\bA^{(N)}) = \frac{1}{\|\euZ\|}\left(\sum_{n=1}^N \|\bG^{(n)}\|^2\right)^{1/2},
\end{align}
where $\bG^{(n)}$ is the mode-$n$ partial derivative of $f$ as defined in \eqref{eq:gradEqn}, we iterate until
\begin{align}\label{eq:stopcrit}
g(\bA^{(1)},\ldots,\bA^{(N)}) < \tau,
\end{align}
or until the maximum number of iterations are reached. In this paper we perform at most 500 iterations of our multilevel method, and at most $10^4$ iterations of ALS. The stopping tolerance is fixed at $\tau = 10^{-10}$.  Table \ref{tab:params} lists the parameters used by the setup and solve phases. As in \cite{Kushnir2010,KahlThesis,SVDAMG}, a larger number of relaxations is required in the setup cycles to produce sufficiently accurate transfer operators.
\begin{table}[!h]
\footnotesize
\caption{CP-AMG-mult and CP-FAS parameters}
\label{tab:params}
\centering
\begin{tabular}{lcc}
\hline\noalign{\smallskip}
Parameter & CP-AMG-mult & CP-FAS\\
\noalign{\smallskip}\hline\noalign{\smallskip}
Pre-relaxations $\nu_1$  & 5 & 1\\
Post-relaxations $\nu_2$ & 5 & 1\\
Relaxations on coarsest level $\nu_c$ & 100 & 50 \\
Cycle type & V-cycle & V-cycle\\
Number of test blocks $n_t$ & 2 & n/a\\
\noalign{\smallskip}\hline
\end{tabular}
\end{table}

For each numerical test we perform ten runs with a different random initial guess for each run. The values reported in the tables represent averages over the successful runs, where a run is deemed successful if the stopping criterion is satisfied prior to reaching the iteration limit. 
The tables compare the ALS method and the multilevel method with or without FMG-CP-FAS as part of the setup phase (see \S\ref{sec:implementation}). For ALS we report the average number of iterations, the average execution time and the number of successful runs. For the multilevel method we report the average number of iterations (setup and solve phases), the average total execution time, the number of successful runs, the average speedup over ALS and the number of levels. The average speedup is determined as follows. For a given test and run, if both ALS and the multilevel method were successful, we divide the execution time of ALS by the execution time of the multilevel method to obtain the speedup for that run. These values are then averaged to obtain the average speedup for that test. We note that execution times do not include the evaluation of the stopping criterion.

\subsection{Implementation details}
\label{sec:implementation}

The multilevel setup and solve phases have thus far been described separately, however, these phases can be combined in the following simple way. Since the factor matrices lie only approximately in the range of the interpolation operators, convergence of the setup cycles, as measured by the functional $g$, should stagnate after a few iterations. Therefore, after each setup cycle the current iterate $\{\bA_{new}\}$ is compared to the previous iterate $\{\bA_{old}\}$ and the setup cycles are halted once
\begin{align}\label{eq:stagnation}
g(\{\bA_{new}\}) > (1-\varepsilon)g(\{\bA_{old}\}),
\end{align}
where the tolerance is set at $\varepsilon = 0.1$. Moreover, at most five setup cycles are performed, and \eqref{eq:stagnation} is checked only after three setup cycles have completed. Once the setup phase has completed, solution cycles are performed until the stopping criterion \eqref{eq:stopcrit} is satisfied. To improve robustness, we can also try to detect stagnation of the solution cycles. After five solution cycles have completed, we check if $g(\{\bA_{new}\}) \geq g(\{\bA_{old}\})$ in each subsequent iteration. If this condition is satisfied then the current iterate $\{\bA_{new}\}$ is discarded and the transfer operators are rebuilt by one down-sweep of CP-AMG-mult with the previous iterate $\{\bA_{old}\}$ used for the boot factors. This process is carried out at most once, and any further indications of stagnation are ignored. Note that the boot factors are not updated by the down-sweep of CP-AMG-mult, as doing so would ruin any progress made by the solution cycles.

The combination of the setup and solve phases described above can be modified to include an FMG-CP-FAS cycle as part of the setup phase. After the setup cycles have completed, we perform one FMG-CP-FAS cycle to compute a new approximation to the boot factors. The transfer operators are then rebuilt using one down-sweep of CP-AMG-mult. Note that while the TFMs are updated by this process, the boot factors are not. We refer to this combination as `Multilevel + FMG' in the tables and figures.

We conclude this section by considering the computational costs of one setup cycle, one solution cycle, and one FMG cycle. Let $\ell = 0,\ldots,L-1$ index the levels, and define
\[
I_n^\ell \coloneqq \frac{I_n}{2^\ell}, \qquad P^\ell \coloneqq \prod_{n=1}^N I_n^\ell = \frac{1}{2^{N\ell}}\prod_{n=1}^N I_n, \qquad S^\ell \coloneqq \sum_{n=1}^N I_n^\ell = \frac{1}{2^\ell}\sum_{n=1}^N I_n
\]
for any $\ell \geq 0$. Assume for simplicity that $\euZ$ is dense and that each mode is coarsened at the same rate with $L$ being the same for each mode.
Consideration of Algorithm \ref{alg:multmlevel} shows that the most expensive operations on each level are the construction of the coarse-level operator, the relaxations, and the construction of the interpolation operators, in particular computing the weights for the least-squares fits.
The coarse-level tensor is constructed by sequentially taking the $n$-mode product of the current tensor with the $n$th restriction operator for $n = 1,\ldots,N$. Computing $\hat{\euZ}^c$ on level $\ell$  requires $\mO(P^\ell S^\ell)$ operations.
The dominant computation for the relaxations and least-squares weights is the matrix product $\bZ_{(n)}\bPhi^{(n)}$. Since $\bZ_{(n)}$ is of size $I_n^\ell \times (P^\ell/I_n^\ell)$ and $\bPhi^{(n)}$ is of size $(P^\ell/I_n^\ell) \times R$ on the $\ell$th level, forming this product requires $\mO(NP^\ell R)$ operations. Therefore, by summing over all the levels, to leading order one setup cycle requires approximately
\begin{align}\label{eq:setupCost}
\left[\left(\frac{2^N}{2^N-1}\right)(n_t(\nu_1+1)+\nu_1+\nu_2+1)+\frac{\nu_c}{2^{N(L-1)}}\right]\cdot\mO(NPR)\nonumber\\
+ \left(\frac{2^N}{2^N-1}\right)\cdot\mO(PS)
\end{align}
operations, where $P = P^0$ and $S = S^0$. We note that $PS$ scales only slightly worse than linear in $P$, and in particular $NPI_{min} \leq PS \leq NPI_{max}$ where $I_{min}$ and $I_{max}$ are the sizes of the smallest and largest modes, respectively. Consideration of Algorithm \ref{alg:CPFAS} shows that the most expensive operations on each level are the relaxations and the construction of the right-hand sides. By a similar analysis, to leading order one solution cycle requires approximately
\begin{align}\label{eq:solnCost}
\left[\left(\frac{2^N}{2^N-1}\right)(\nu_1+\nu_2+1+1/2^N)+\frac{\nu_c}{2^{N(L-1)}}\right]\cdot\mO(NPR)
\end{align}
operations. Similarly, it follows that to leading order one FMG-CP-FAS cycle requires approximately
\begin{align}\label{eq:FMGCost}
\left[\left(\frac{2^N}{2^N-1}\right)^2(\nu_1+\nu_2+1+1/2^N)+\frac{\nu + (L-1)\nu_c}{2^{N(L-1)}}\right]\cdot\mO(NPR)
\end{align}
operations. Note that in \eqref{eq:FMGCost} $\nu$ is the number of ALS iterations performed on the coarsest level and $(\nu_1,\,\nu_2,\,\nu_c)$ are the CP-FAS parameters.
In general a solution cycle is significantly cheaper than a setup cycle because of the extra work required by a setup cycle to relax on the TFMs, the typically larger number of relaxations performed on each level of a setup cycle, and the added work of constructing the coarse-level tensors (i.e., the $\mO(PS)$ term).
If $\euZ$ is sparse then further savings are possible on the finest level. In particular, to leading order the cost of one relaxation reduces to $NR$ times the number of nonzero elements in $\euZ$.
In our current framework the coarse tensors will in general be dense; multiplication by the inverted Cholesky factors as in \eqref{eq:ctensor_transformed} destroys any sparsity. Therefore, it may be interesting to consider alternative formulations of the coarse-level equations, for example, working directly with equations of the form in \eqref{eq:coarseeq1}; see also \cite{SVDAMG}.

\subsection{Sparse tensor test problem}
\label{sec:spten}

The first test problem we consider is the standard finite difference Laplacian tensor on a uniform grid of size $s^d$ in $d$ dimensions. This test problem yields an $N$-mode sparse tensor $\euZ$ of size $s\times s\times\cdots\times s$ with $N = 2d$. We can efficiently construct $\euZ$ by reshaping the $s^d \times s^d$ matrix
\[
\bZ = \sum_{k=1}^d \bI_{\ell(k)}\otimes\bD\otimes \bI_{r(k)},
\]
where $\bI_{r(k)}$ is the $s^{k-1} \times s^{k-1}$ identity matrix, $\bI_{\ell(k)}$ is the $s^{d-k} \times s^{d-k}$ identity matrix, and $\bD$ is the $s \times s$ tridiagonal matrix with stencil $[-1,\,2,\,-1]$. While this test problem is somewhat pedagogical in nature, it offers a good starting point to illustrate our method. The parameters for the various test tensors that we consider are given in Table \ref{tab:LapParams}.

\begin{table}[!h]
\footnotesize
\caption{Parameters for sparse problem.}
\label{tab:LapParams}
\centering
\begin{tabular}{rc}
\hline\noalign{\smallskip}
test & problem parameters\\
\noalign{\smallskip}\hline\noalign{\smallskip}
 1 & $N = 4,\,s = 20,\,R = 4$\\
 2 & $N = 4,\,s = 20,\,R = 5$\\
 3 & $N = 4,\,s = 20,\,R = 6$\\
 4 & $N = 4,\,s = 50,\,R = 2$\\
 5 & $N = 4,\,s = 50,\,R = 3$\\
 6 & $N = 4,\,s = 50,\,R = 4$\\ 
 7 & $N = 6,\,s = 20,\,R = 2$\\
 8 & $N = 6,\,s = 20,\,R = 3$\\
 9 & $N = 6,\,s = 20,\,R = 4$\\
10 & $N = 8,\,s = 10,\,R = 2$\\
11 & $N = 8,\,s = 10,\,R = 3$\\
\noalign{\smallskip}\hline
\end{tabular}
\end{table}

\begin{table}[!h]
\footnotesize
\caption{Sparse problem. Average number of iterations and time (in seconds) until the stopping criterion is satisfied with stopping tolerance $10^{-10}$. Here `it' is the number of iterations, `spd' is the multilevel speedup compared to ALS, `ns' is the number of successful runs, and `levs' is the number of levels.}
\label{tab:Laplacian}
\centering
\begin{tabular}{rrrrrrrrrrrrr}
\hline\noalign{\smallskip}
 &  \multicolumn{3}{c}{ALS} & \multicolumn{4}{c}{Multilevel} & \multicolumn{4}{c}{Multilevel + FMG} &\\
\noalign{\smallskip}\hline\noalign{\smallskip}
test & it & time & ns & it & time & spd & ns & it & time & spd & ns & levs\\
\noalign{\smallskip}\hline\noalign{\smallskip}
 1 & 1897 &  24.0 & 10 &  37 &   7.6 &  3.2 & 10 &   36 &   8.0 &  3.1 & 10   &   2\\
 2 & 3329 &  53.6 &  8 &  64 &  15.1 &  4.5 & 10 &   42 &  11.8 &  5.1 & 10   &   2\\
 3 & 3587 &  70.3 &  9 &  67 &  17.5 &  4.0 &  9 &   32 &  10.6 &  6.8 & 10   &   2\\
 \noalign{\smallskip}\hline\noalign{\smallskip}
 4 & 5457 & 105.9 & 10 & 123 &  40.5 &  2.7 & 10 &  120 &  41.2 &  2.6 & 10   &   4\\
 5 & 5508 & 150.3 &  4 & 182 &  64.9 &  2.3 &  9 &   99 &  40.9 &  3.8 &  9   &   4\\
 6 & 6788 & 244.0 &  3 & 150 &  62.0 &  5.2 & 10 &  136 &  58.7 &  5.4 &  9   &   4\\
\noalign{\smallskip}\hline\noalign{\smallskip} 
 7 & 1619 & 187.0 & 10 &  48 & 111.4 &  1.7 & 10 &   52 & 128.3 &  1.5 & 10   &   3\\
 8 & 3481 & 610.0 & 10 &  72 & 164.8 &  3.7 & 10 &   70 & 177.9 &  3.5 & 10   &   3\\
 9 & 4085 & 939.5 & 10 &  76 & 209.6 &  4.5 & 10 &   78 & 228.9 &  4.2 & 10   &   3\\
\noalign{\smallskip}\hline\noalign{\smallskip}
10 &  634 & 229.5 & 10 &  48 & 170.8 &  1.3 & 10 &   56 & 203.5 &  1.1 & 10   &   3\\
11 & 1743 & 943.3 & 10 &  39 & 402.3 &  2.3 & 10 &   43 & 474.0 &  2.0 & 10   &   3\\
\noalign{\smallskip}\hline
\end{tabular}
\end{table}

The results of the tests for the sparse test case are given in Table \ref{tab:Laplacian}. The results show that our multilevel approach is anywhere from two to seven times faster than ALS for this test problem. For tests 1 to 6 (order 4 tensors), larger speedups are observed for the multilevel method with FMG. However, for tests 7 to 11 (order 6 and 8 tensors) larger speedups are observed for the multilevel method without FMG. For higher-order tensors the setup phase of multilevel method with FMG is considerably more expensive than the setup phase of the multilevel method without FMG. The multilevel variants demonstrate similar robustness to varying initial guesses for this problem, however, in general we expect the multilevel method with FMG to be the most robust option. We also observe the trend that for each grouping of tests in Table \ref{tab:Laplacian}, the speedup tends to increase as the number of components $R$ increases.

\begin{figure}[!h]
\centering
\includegraphics{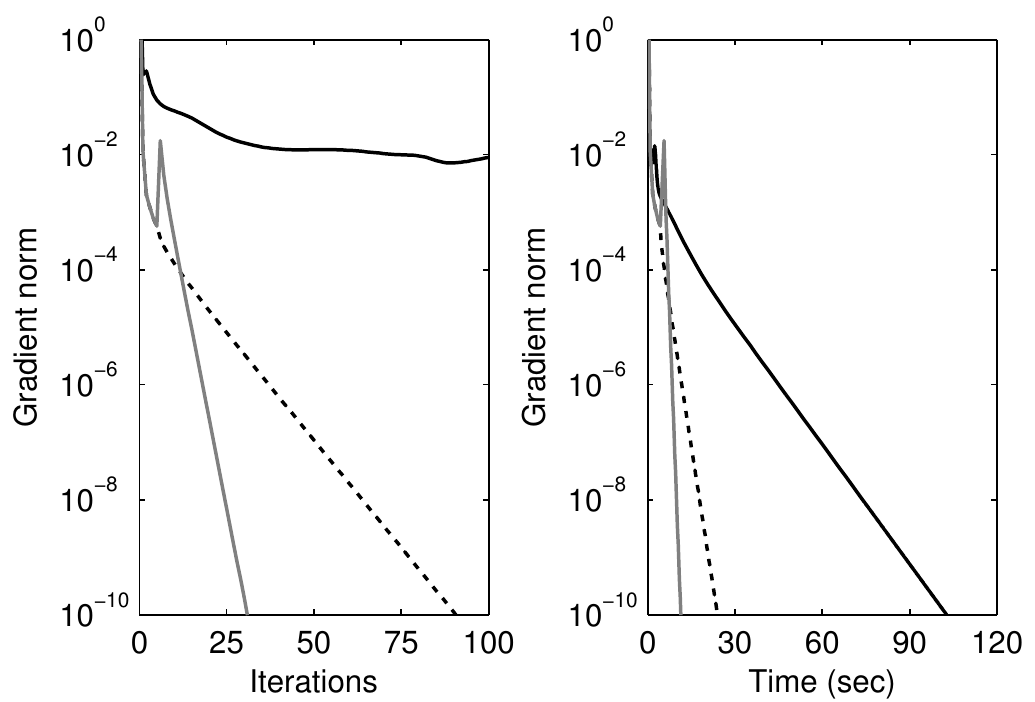}
\caption{Sparse problem. Convergence plot for test 3 from Table \ref{tab:LapParams} ($N = 4$, $s = 20$, $R = 6$). The solid black line is ALS, the solid gray line is the multilevel method with FMG, and the dashed line is the multilevel method without FMG.}
\label{fig:Lap3}
\end{figure}

\begin{figure}[!h]
\centering
\includegraphics{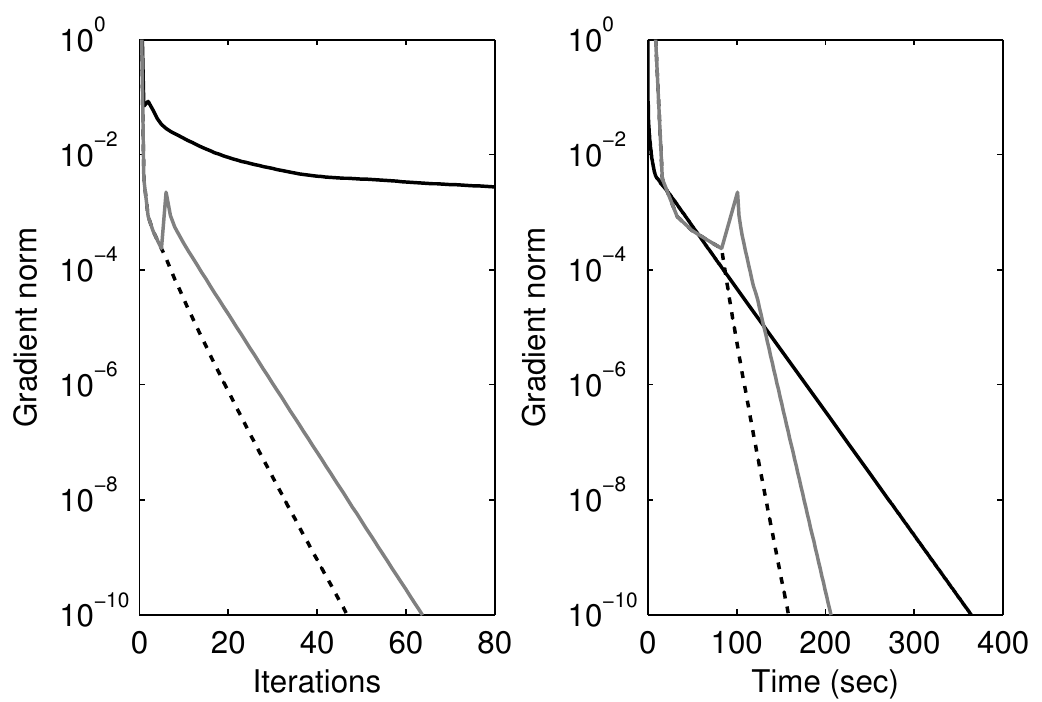}
\caption{Sparse problem. Convergence plot for test 9 from Table \ref{tab:LapParams} ($N = 6$, $s = 20$, $R = 4$). The solid black line is ALS, the solid gray line is the multilevel method with FMG, and the dashed line is the multilevel method without FMG.}
\label{fig:Lap9}
\end{figure}

Figures \ref{fig:Lap3} and \ref{fig:Lap9} illustrate the convergence history of ALS and the multilevel method for one run of tests 3 and 9 in Table \ref{tab:Laplacian}, respectively. These plots are typical of the performance observed for this test problem. We note that the spike in the `Multilevel + FMG' curves is due to the initial approximation to the solution computed by the single FMG-CP-FAS cycle performed after the setup cycles.


\subsection{Dense tensor test problem}
\label{sec:denseten}

The second test problem we consider is a dense, symmetric third-order tensor $\euZ \in \mathbb{R}^{s\times s\times s}$ whose elements are given by
\[
z_{ijk} = \big(i^2+j^2+k^2\big)^{-1/2}\quad\text{for } i,j,k = 1,\ldots,s.
\]
This tensor was used as a test case in \cite{Oseledets2006}, which compares various methods for computing the CP decomposition including ALS. It was also considered in \cite{Oseledets2008}, which describes a novel method for computing the Tucker decomposition of third-order tensors. As mentioned in \cite{Oseledets2008}, $\euZ$ arises from the numerical approximation of an integral equation with kernel $1/\|\bx-\by\|$ acting on the unit cube and discretized by the Nystr\"{o}m method on a uniform grid.
In this section we compute CP decompositions for $R = 2,\,3,\,4,\,5$. It has been observed numerically that when $R \geq 4$, ALS may be extremely slow to converge, for some initial guesses requiring on the order of $10^5$ iterations, with highly non-monotonic convergence behavior. The performance of our method when $R \geq 4$ is less robust than desired because the multigrid framework uses a single interpolation operator for each factor matrix. Even so, depending on the initial guess our method may still demonstrate a significant improvement over ALS.

\begin{table}[!h]
\footnotesize
\caption{Dense problem. Average number of iterations and time (in seconds) until the stopping criterion is satisfied with stopping tolerance $10^{-10}$. Here `it' is the number of iterations, `spd' is the multilevel speedup compared to ALS, `ns' is the number of successful runs, and `levs' is the number of levels.}
\label{tab:IJK}
\centering
\begin{tabular}{rcrrrrrrrr}
\hline\noalign{\smallskip}
 &  & \multicolumn{3}{c}{ALS} & \multicolumn{4}{c}{Multilevel + FMG} &\\
\noalign{\smallskip}\hline\noalign{\smallskip}
test & \multicolumn{1}{c}{problem parameters} & it & time & ns & it & time & spd & ns & levs\\
\noalign{\smallskip}\hline\noalign{\smallskip}
 1 & $s = \phantom{1}50,\,R = 2$ &  161 &    0.7 & 10  &    7 &  2.2 &   0.3 & 10  &    5\\
 2 & $s = \phantom{1}50,\,R = 3$ & 2435 &   11.8 & 10  &   15 &  2.7 &   4.7 & 10  &    5\\
 3 & $s = \phantom{1}50,\,R = 4$ & 4838 &   26.1 &  5  &   97 & 10.1 &   4.4 &  7  &    4\\
 4 & $s = \phantom{1}50,\,R = 5$ &      &        &     &  301 & 30.0 &       &  3  &    4\\
\noalign{\smallskip}\hline\noalign{\smallskip}
 5 & $s = 100,\,R = 2$ &  253 &   10.7 & 10  &    7 &  7.9 &   1.4 & 10  &    6\\
 6 & $s = 100,\,R = 3$ & 1695 &   80.2 &  9  &    9 &  8.1 &  10.1 & 10  &    6\\
 7 & $s = 100,\,R = 4$ & 3836 &  202.2 &  6  &   82 & 27.5 &  14.1 &  9  &    5\\
 8 & $s = 100,\,R = 5$ & 7854 &  455.2 &  2  &  192 & 62.0 &   9.0 &  4  &    5\\
\noalign{\smallskip}\hline\noalign{\smallskip} 
 9 & $s = 200,\,R = 2$ &   274 &   90.3 & 10  &    7 &  48.8 &   1.9 & 10  &    7\\
10 & $s = 200,\,R = 3$ &  1830 &  682.3 & 10  &   12 &  61.7 &  11.2 & 10  &    7\\
11 & $s = 200,\,R = 4$ &  2998 & 1249.5 &  8  &   79 & 178.0 &  11.6 &  9  &    6\\
12 & $s = 200,\,R = 5$ &  5686 & 2611.4 &  3  &  220 & 440.7 &   2.6 &  4  &    6\\
\noalign{\smallskip}\hline
\end{tabular}
\end{table}

The results of the tests for the dense test case are given in Table \ref{tab:IJK}. We note that only the multilevel method with FMG is considered for this test problem (see the description in \S\ref{sec:implementation}). The blank entries for test 4 in Table \ref{tab:IJK} indicate that ALS did not have any successful runs. 
For $R \geq 3$ our multilevel approach can lead to significant savings in iterations and execution time. The speedup is less impressive when $R = 2$, since ALS already converges quickly without any multigrid acceleration. It is also apparent that as the number of components increases, the number of successful runs of the multilevel method, and of ALS, decreases. For initial guesses in which the multilevel method failed to converge, there was typically a rapid decrease in the gradient norm, followed by convergence stagnation of the solution cycles. This behavior suggests that the setup phase was unable to construct transfer operators that adequately represented the solution in their range. Such cases were also characterized by slow convergence of ALS.
\begin{figure}[!h]
\centering
\includegraphics{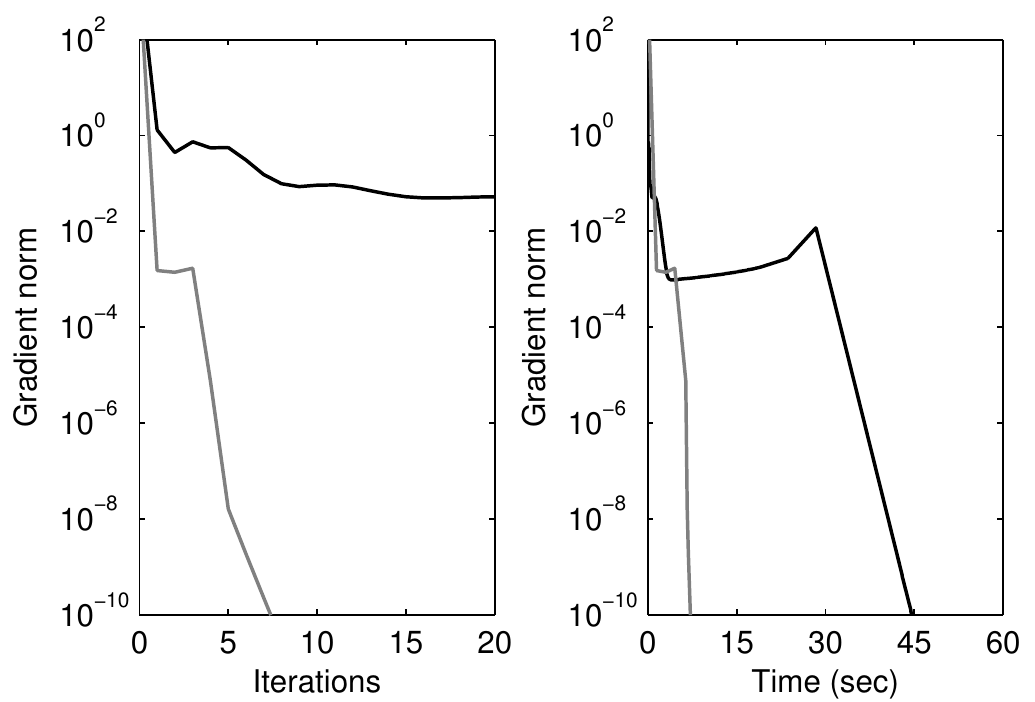}
\caption{Dense problem. Convergence plot for test 6 from Table \ref{tab:IJK} ($s = 100$, $R = 3$). The solid black line is ALS and the solid gray line is the multilevel method with FMG.}
\label{fig:IJK6}
\end{figure}
\begin{figure}[!h]
\centering
\includegraphics{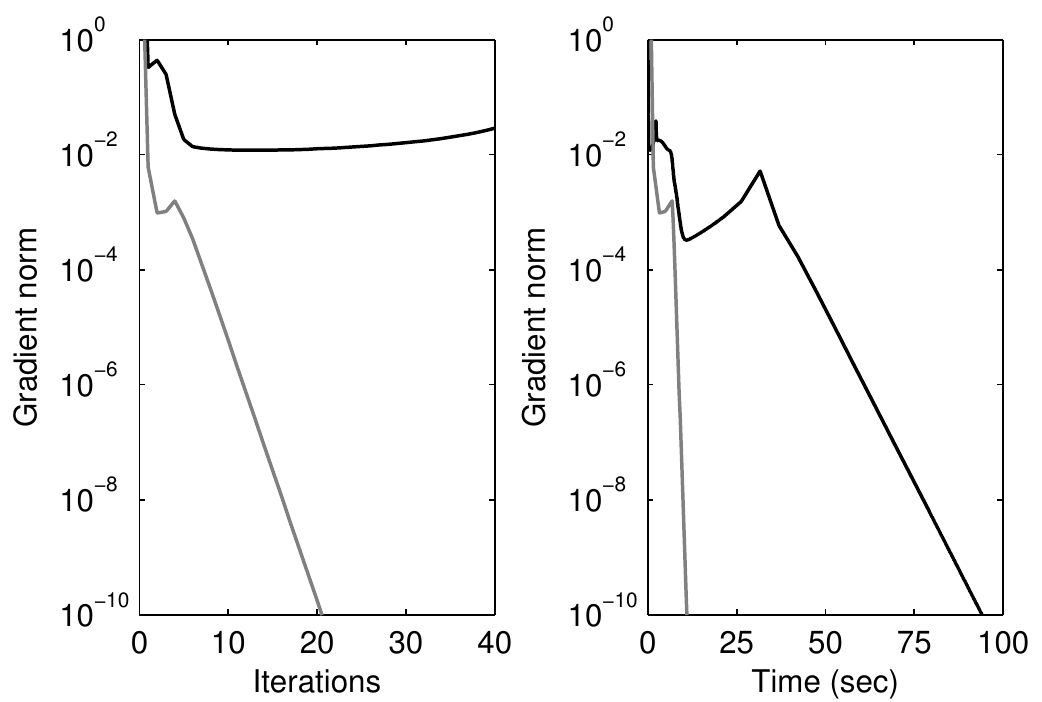}
\caption{Dense problem. Convergence plot for test 7 from Table \ref{tab:IJK} ($s = 100$, $R = 4$). The solid black line is ALS and the solid gray line is the multilevel method with FMG.}
\label{fig:IJK7}
\end{figure}
\begin{figure}[!h]
\centering
\includegraphics{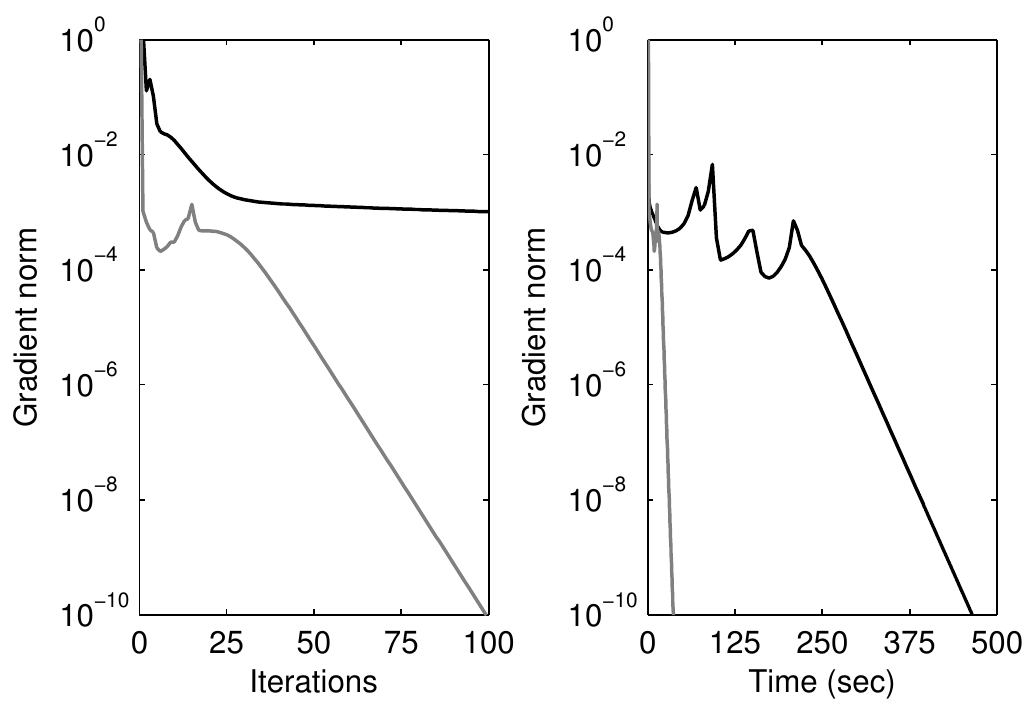}
\caption{Dense problem. Convergence plot for test 8 from Table \ref{tab:IJK} ($s = 100$, $R = 5$). The solid black line is ALS and the solid gray line is the multilevel method with FMG.}
\label{fig:IJK8}
\end{figure}

Figures \ref{fig:IJK6}, \ref{fig:IJK7} and \ref{fig:IJK8} illustrate the convergence history of ALS and the multilevel method for one run of tests 6, 7 and 8 in Table \ref{tab:IJK}, respectively. Figure \ref{fig:IJK8} $(R = 5)$ shows how  ALS can initially be slow to converge with erratic convergence behavior: for the first half of the run its gradient norm fluctuates with little decrease. Such behavior can make it very difficult for the setup phase to construct adequate transfer operators.

\section{Concluding remarks}
\label{sec:conc}

We have presented a new algorithm for computing the rank-$R$ canonical decomposition of a tensor for small $R$. As far as we are aware our method is the first genuine multigrid algorithm for computing the CP decomposition. Our work is also significant in that it presents the first adaptive AMG method for a nonlinear optimization problem. Similar to the method presented in \cite{SVDAMG} for computing SVD triplets of a matrix, we combined an adaptive multiplicative setup phase with an additive solve phase. The ALS method was used as the relaxation scheme. Numerical tests with dense and sparse tensors of varying sizes and orders (up to order 8) that are related to PDE problems showed how our multilevel method can lead to significant speedup over standalone ALS when high accuracy is desired.

Avenues of further research are plentiful. For example, it may be worthwhile to investigate a more sophisticated setup phase that iteratively combines the CP-AMG-mult cycles and CP-FAS-FMG cycles. As discussed briefly in \S\ref{sec:implementation}, an alternative formulation of the coarse-level equations without the inverted Cholesky factors may be fruitful for sparse problems. In addition to PDE-related tensors, there may be other classes of tensors for which multigrid acceleration of ALS may be beneficial, but identifying and studying such classes remains a topic of future research. It would also be interesting to consider other ALS-type methods for the relaxations, for example, those using line searches, as well as to apply our method to the regularized optimization formulation of CP as described in \cite{CPOPT}. Similarly, it would be interesting to generalize our multilevel framework to other similar tensor optimization problems such as the Tucker decomposition \cite{TensorReview}, block tensor decompositions \cite{Lathauwer2008II,Lathauwer2008III}, best rank-$(R_1,\ldots,R_N)$ approximations \cite{Lathauwer2000,Ishteva}, and to other nonlinear optimization problems.
\newpage

\bibliographystyle{siam}
\bibliography{MGtensorCP}
\end{document}